# Women's participation in mathematics in Scotland, 1730–1850

Amie Morrison & Isobel Falconer

Published online: 25 Apr 2022.

Submit your article to this journal

View related articles

View Crossmark data

# Women's participation in mathematics in Scotland, 1730–1850

Amie Morrison
*University of St Andrews, St Andrews, Fife, Scotland**

Isobel Falconer
*School of Mathematics and Statistics, University of St Andrews, St Andrews, Fife, Scotland*

The eighteenth century saw a flourishing of scientific and philosophical thought throughout Scotland, known as the Scottish Enlightenment. The accomplishments of prominent male figures of this period have been well documented in all disciplines. However, studies of women's experiences are relatively sparse. This paper partially corrects this oversight by drawing together evidence for women's participation in mathematics in Scotland between 1730 and 1850. In considering women across all social classes, it argues for a broad definition of 'mathematics' that includes arithmetic and astronomy, and assesses women's opportunities for engagement under three headings: education, family, and sociability. It concludes that certain elements of Scottish Enlightenment culture promoted wider participation by women in mathematical activities than has previously been recognized, but that such participation continued to be circumscribed by societal views of the role of women within family formation.

## Introduction

The eighteenth century bore witness to a flourishing of scientific and philosophical thought throughout Scotland, known as the Scottish Enlightenment. The accomplishments of prominent male figures of this period have been well-documented in all disciplines. Unfortunately, it remains, as Rosalind Russell notes, that 'Scotswomen are underrepresented', and studies of women's experiences are relatively sparse (Russell 1988, vii). Moreover, since, as Russell details, women were actively discouraged from scientific endeavour by eighteenth-century Scottish society, and could not enter universities,[1] existing scholarship has neglected the study of Scottish women in mathematics.[2] Aside from biographical accounts of extraordinary individuals like Mary Somerville, we have found no attempt to survey the extent and experiences of Scottish women obtaining a mathematical education, practicing mathematics, or contributing to mathematical research. This paper will attempt to correct this oversight, assessing the evidence for women's participation in mathematics in Scotland between 1730 and 1850.

---

*School of Mathematics and Statistics, University of St Andrews, St Andrews, Fife, Scotland. Email: amie.morrison@live.com

[1] All four Scottish Universities began admitting women to degree courses for the first time in 1892, following the provisions of the 1889 Universities (Scotland) Bill. Beginning in 1867, limited access to higher education for women was provided by the Edinburgh Ladies Educational Association (founded 1867), the similar Aberdeen Ladies Educational Association, the LLA (Lady Literate in Arts) distance learning scheme started by St Andrews in 1877, some specific Ladies Classes at Glasgow (began 1868), and two Colleges: Queen Margaret (women only) in Glasgow (1883), and Dundee University College (1881).

To do so, it becomes necessary to be flexible in defining 'mathematics' as a field. Although this study will discuss standard areas of mathematics like algebra and geometry, a fuller picture of women's mathematical engagement can be obtained by allowing the definition of 'mathematics' to include elementary arithmetic and the numerical aspects of astronomy; most evidence of activity is at these levels and there are few indications of engagement at higher levels during the period, perhaps due to the lack of access to higher education mentioned above. The period covered by this paper will extend from 1730 to 1850, incorporating the changes to Scottish society wrought by the Enlightenment. Rosalind Carr defines the Scottish Enlightenment as an 'often disparate ideological and cultural movement unified by a discourse of improvement', where 'improvement' can be understood as 'an imperative to achieve and maintain social progress' (Carr 2014, 1). By extending the period of study beyond the traditional periodization of the Scottish Enlightenment to include the first half of the nineteenth century, it will be possible to encompass a wide range of women's experiences, and to consider whether this 'imperative' for social change had any impact on societal attitudes towards women in mathematics.

The findings of this paper are presented under three headings: Education, Family and Sociability. The first section will discuss the forms of mathematical education available to women of all social classes during this period. Through surveying school advertisements in Scottish newspaper archives, we have identified many women's schools that offered classes in arithmetic, and a few exceptional establishments that offered separate classes in mathematics. Other forms of education, such as public lectures and schools for adult women, are also incorporated, from various other primary and secondary sources. The section on Family examines how women could be put into a position where arithmetical, and occasionally higher mathematical, skills were required of them within the family, particularly through the task of accounting in the family home or business. It draws on examples from a substantial piece of research carried out by Rosalind Russell on the lives of business-owning Scotswomen during this period, anecdotes from a multitude of other secondary sources, and analysis of an arithmetic textbook written for women in Scotland. It also discusses continuity and change in the concept of the family during this period and compares similar situations occurring elsewhere in Britain and Europe. Building upon the work of Ruth Watts, Mary Brück and Brigitte Stenhouse, the final section on Sociability pieces together several under-explored primary and secondary sources to investigate the conversational culture of the Enlightenment, and how this can be used to understand the processes through which educated women engaged with mathematics in Scotland.

Through examining how women could acquire a mathematical education, practice mathematics through family life, and participate in mathematical discussion through social connections, this paper will demonstrate that women in eighteenth- and nineteenth-century Scotland participated in mathematical activity far more than might previously have been assumed. Increasingly, though, they did so in a role as 'help-meet' to a male family head, and we will maintain throughout that Enlightenment discourse did little to change societal attitudes towards women pursuing mathematical study; greater engagement with mathematics throughout this period owed more to favourable social and political circumstances than it did to enlightened educational reform.

---

[2]An exception, but one that focuses primarily on *girls* in the 1850–1870 period, immediately after ours, is (Moore 1984).

**Education**

During the early eighteenth century, opportunities to learn mathematics through formal educational routes were few, even at a rudimentary level. Many young women of the lower social classes were illiterate and received no education whatsoever, whereas upper-class families might offer private tuition to their daughters in 'accomplishments' such as musical instruments, dancing, or languages (Hughes 2014). However, as the eighteenth century progressed, the number of private girls' schools increased rapidly (Moore 2003). Analysis of archived newspaper advertisements illustrates that several of these offered mathematical tuition for young women. As early as 1786, an advertisement mentions that 'Mrs Barker intends opening a Day School for the Education of Young Ladies' in Edinburgh, where 'every suitable branch of Education' would be taught, including arithmetic (Education for Young Ladies', 1786). Another school in Edinburgh run by 'Mrs Crewe and Daughters from London' in 1800 offered arithmetic classes, in addition to other subjects, for forty guineas per annum (the equivalent of around £1,850 today) (English Ladies School', 1800). Evidently, obtaining formal education in arithmetic or any other subject for young women at the turn of the nineteenth century did not come at a low cost, and would therefore primarily be available to upper-class families. What is not clear from these advertisements is what the term 'arithmetic' implied; as at the parish schools discussed below, the level may possibly have been surprisingly high in some cases.

Private girls' schools proliferated during the early nineteenth century, particularly in the 1830s. Arithmetic was offered for boarders at 'Mrs Kiloh's Ladies' School' in Aberdeen and for pupils at Gayfield Square Ladies' School in Edinburgh (Ladies School, 1836; Gayfield Square Ladies' School, 1835). Several private tutor groups were also advertised, like Mr Mcfarlane and Mr Gouinlock's classes in 'Writing, Arithmetic and Book-keeping', offered in 1828 (Writing, Arithmetic, and Book-keeping, 1828). These classes even offered a Dux prize for the best two girls in arithmetic, a practice which was also taken up by the Inverness Royal Academy, Perth Grammar School and Leith High School during this period (Inverness Royal Academy, 1839; Perth Seminaries, 1846; Leith High Schools, 1833). Katharine Glover notes that upper-class families often educated young men at boarding schools in London and this was 'increasingly extended to Scottish girls in the late eighteenth century' (Glover 2010, cited in Whiting 2018, 26). While their brothers studied a liberal arts curriculum which, if it included mathematics, comprised mainly Euclid and algebra, the girls, as in Scotland, were generally limited to arithmetic and accounting. However, several progressive London schools went further; Margaret Bryan's school in Blackheath was a notable example. Bryan included subjects such as trigonometry, optics, gravity, mechanics, pneumatics, and hydrostatics in her lessons as well as astronomy and the use of the globes (Bryan, 1799; Benjamin 1991, 38; Keene 2011, 533).

There were also opportunities available for women of lower social standing. 'Parish schools' in rural areas and public 'burgh schools' in towns were, in principle if not always in practice, open to girls as well as boys. By 1851, over one-third of the children learning arithmetic in parish schools were girls and at least forty-five elementary public schools were teaching mathematics to girls (Moore 1984, 125–129). Even the 'arithmetic' classes might encompass Euclid, mensuration, navigation and surveying in some of these schools (Gray 1952, 37). By 1826, at least 140 parish and burgh schools, mainly around the coast, offered navigation; an outstanding example was at Lecropt, near Dunblane, where fluxions, navigation, astronomy and plane and spherical trigonometry were apparently taught, although it is not clear that any girls took these options (Gavine, 1990, 7). Throughout most of the period, girls did not have access to the other avenues for learning practical mathematics that were available to boys, such as the town academies and navigation schools detailed by Gavine (1990). However, by the mid nineteenth century, things were changing: in 1848, the Public Commercial and Mathematical School in Aberdeen offered classes in arithmetic and bookkeeping for young ladies (Summer Classes, 1848). A record of a council board meeting in 1847 accounts for a decision to grant a £25 salary to a ladies' seminary in Forres, which included instruction in arithmetic, and was aimed at 'the daughters of respectable tradesmen and merchants, at such moderate fees as they could afford to pay' (Ladies' School, 1847). Houses or Schools of Industry set up by charitable societies, like the Leith Female Charity School of Industry in 1802, allowed young poor women to acquire industrial skills like spinning or sewing, and were often run by committees of upper-class women (Rendall 2017, 213). Such schools occasionally offered arithmetic classes, such as another girls' school in Leith in 1819 where women were given classes in writing and arithmetic as 'rewards of good behaviour', in addition to 'reading, sewing, knitting and spinning' (ibid.). Although this school might appear progressive in offering such subjects, it demonstrates that arithmetical skills were ultimately considered superfluous for women of these social classes, being offered as a 'reward'. Nonetheless, the presence of formal educational establishments where young women of all social classes could learn mathematics at some level had multiplied by 1850, undoubtedly encouraging a higher rate of mathematical literacy.

One establishment notable for its mathematical instruction was the Scottish Institution for the Education of Young Ladies (Macdonald and Hope 1993; Moore 2003, 251). Founded in 1834 and lasting until 1871, the Scottish Institution was a response to the educational reform movement relating to science sweeping Scotland at the time. The school offered classes not just in arithmetic (taught by a 'Mr Trotter'),[3] but also a separate class in mathematics, taught by George Lees, who also taught at the Scottish Military and Naval Academy (Scottish Institution for the Education of Young Ladies, 1835). Donald Sutherland notes that the syllabus of Mr Lees' lectures on mathematics included 'Properties of Matter, Elements of Mechanics, Pendulums, Central Forces,

---

[3] We are grateful to a referee for the suggestion that 'Mr Trotter' may refer to Alexander Trotter, son of James Trotter. Alexander evidently wrote many of the exercises in his father's *Lessons in Arithmetic for Junior Classes* and was the author of *A Key to [J.] Trotter's Complete System of Arithmetic*. Both father and son were connected with the Scottish Naval and Military Academy.

Movements of Heavenly Bodies, Motion and Pressure of Air, Acoustics, Theory of Music', all of which were relatively advanced topics considering the novelty of mathematical instruction for the female attendees (Sutherland 1938, 176). Lees held progressive views on education, and published a textbook on mathematics for students at the Edinburgh School of Arts, an institution aimed at providing education to adults in scientific trades (Kelly 1952, 24–25; Lees 1826). As will be discussed shortly, these institutions were often welcoming to women, which might explain Lees' willingness to engage in women's mathematical education through the Scottish Institution. As stated by its founders, the purpose of the Scottish Institution was to create the 'union of the scientific with the ornamental branches of education, and the accessibility to such a liberal course by a moderate fee' (Scottish Institution for the Education of Young Ladies, 1849). Instead of pricing each subject individually, the course fee was fixed at 'five guineas per quarter' for all classes, which might have encouraged women to attend classes that their families thought were unnecessary (Scottish Institution for the Education of Young Ladies, 1835). This may have partially ensured attendance at the mathematical and scientific classes, which would not traditionally have been oriented towards women. Other educational establishments followed suit in Edinburgh, with the 'Edinburgh Ladies' Institution' and 'Mrs Ponsonby's Establishment for the Board and Education of Young Ladies' each offering branches in 'arithmetic' and 'mathematics' (Edinburgh Ladies' Institution, 1846; Mrs Ponsonby's Establishment for the Board and Education of Young Ladies, 1845; Moore 2003, 251).

It is, however, important not to misinterpret the presence of mathematical education in these schools as a reflection of progressive attitudes towards women. In a pragmatic sense, this education allowed women to obtain greater mathematical knowledge, provided a counterexample to the stereotype that they were unsuited for scientific study, and presented an opportunity to demonstrate their capabilities in this area. However, it was not for these purposes that it was introduced. Mathematics was seen as a subject that could enrich a woman's understanding of the world, not for her own sake, but for that of her children and marriage prospects. The Scottish Institution stated that it educated young women for the sole purpose of creating better mothers for society's future children:

> If the mind of the mother were initiated in the study of those subjects which render the works of the universe around her intelligible, she could point out to her children, not merely their external beauty, but that wise and more hidden arrangement, which expands the flower, and variegates the rainbow, and works in a thousand common objects (Scottish Institution for the Education of Young Ladies, 1835).

Katharine Glover argues that the notion that 'girls should be able to provide intelligent company (particularly within the family)' was 'probably the most effective factor in extending young women's access to knowledge during the eighteenth century (Glover 2011, cited in Barclay and Carr 2018, 177). Any educational benefits for the young women in question were 'frequently incidental, if not always insignificant' (ibid.). Although mathematical education provided a novel opportunity for young women and began to normalize the existence of mathematical knowledge amongst women in society, it is important to note the caveat that such education was introduced largely to help them fulfil their predetermined role as daughters, wives, and mothers.

Although universities did not offer mathematical classes for women, or permit women to attend those for men, until the late nineteenth century, there were instances of public lectures open to women. In the late eighteenth century, James Ferguson, an astronomer from Aberdeen, delivered lectures throughout Britain and always welcomed women, noting that it 'gives him great pleasure to find that many Ladies form themselves into companies, and attend his Lectures' (Brück 2009, 12).[4] He is known to have lectured at least five courses in Edinburgh between March and July 1768 (Millburn 1985, 406). Ferguson published *Astronomy Explained Upon Sir Isaac Newton's Principles, and Made Easy to Those Who Have Not Studied Mathematics* in 1756, and *The Young Gentleman and Lady's Astronomy, Familiarly Explained in Ten Dialogues between Neander and Eudosia* in 1768. The latter explicitly identifies women as part of its readership, and was structured around an imagined conversation, a popular style for women's scientific textbooks during this period. According to the author, neither book assumed prior knowledge of mathematics; the *Young Gentleman's and Lady's Astronomy* aimed to show that men and women could 'acquire a competent knowledge of Astronomy, without any previous knowledge of Geometry or Mathematics' (Ferguson 1768, xi). However, the books did assume basic arithmetic and familiarity with terms such as 'diameter' and 'chord', and introduced concepts such as proportion. The elementary level of mathematics does not appear to be due to a particular lack of esteem for women's abilities; rather, Ferguson aimed his books at a beginner level for both men and women, expecting neither to have received an education beyond arithmetic. This point is made explicit in *The Young Gentleman's and Lady's Astronomy*, where Neander states that men 'very absurdly imagine, because they know nothing of science themselves, that it is beyond the reach of women's capacities', disavowing the notion that women were biologically incapable of it (76). Eudosia does indeed carry out a calculation at one point in the book, computing the speed of light:

> you taught me not only the four common rules of arithmetic before you went to university, but even the rule of three. The Sun's distance from the Earth is 95 millions of miles, in round numbers; and light moves through that space in 8 minutes of time; divide, therefore, 95,000,000 by 8, and the quotient is 11,875,000, for the number of miles that light moves in a minute. Now, I remember that you told me, a cannon-ball moves at the rate of 480 miles in an hour, which is 8 miles in a minute; I therefore divide 11,875,000 by 8, and the quotient is 1,484,375; so that light moves more than a million of times as swift as a cannon-ball — Amazing indeed! (82–83)

Ferguson clearly believed that women were capable of quantitative astronomy, and that it was a suitable subject of study for them. His motivations for this belief are made evident in the text when Neander states that scientific pursuits would give women a 'rational way of spending their time at home' instead of 'murdering it, by going abroad to card-tables, balls, and plays' (75). He notes that 'how much better wives, mothers, and mistresses they would be, is obvious to the common sense of mankind' (76). Echoing those sentiments expressed previously by the Scottish Institution, Ferguson desired a scientific education for women to help them better fulfil their feminine roles. Mary Brück argues that he possessed the 'ideal of women as educated companions to their enlightened menfolk in the family circle'; education in astronomy would enhance their marriageability, much like the traditional feminine accomplishments of music, drawing or dancing (Brück 2009, 15). Regardless of his motivations, he clearly believed, contrary to many authors of his time, that men and women were equally capable of understanding the mathematics of astronomy, and that it was an appropriate pastime for them to engage in.

---

[4]For Ferguson's career see (Millburn and King 1988; Rothman 2000).

A similar course of public lectures was held by John Anderson, Professor of Natural Philosophy at the University of Glasgow from 1757 to 1796, for mixed groups of men and women (Anonymous, 'John Anderson'; Wood 2004). Anderson bequeathed his estate to the foundation of 'Anderson's Institution', which would act as an accessible alternative to the University of Glasgow (Anonymous, 'John Anderson'). His will instructed that the Institution should give a course.

> at least once every year, to be called "The Ladies Course of Physical Lectures", in which no mathematical reasoning shall be used […] the Audiences shall consist of both Ladies and Gentlemen (Pike 2012).

Anderson's insistence on the absence of mathematical reasoning from the proposed ladies' course demonstrates that, unlike Ferguson, he doubted women's ability to understand mathematics, although whether due to lack of prior education or lack of intellect is not clear. However, his advocacy of accessible higher education was reflected in the ethos of the Institution, which by the 1840s offered classes to female students specifically in mathematics (ibid.). In 1840, the Natural Philosophy class contained 130 women (ibid.).

A dispute at Anderson's Institution in 1823 led to the creation of an offshoot called the Mechanics' Institution, which became one of numerous 'Mechanics' Institutes', or 'Schools of Arts', aimed at providing a scientific education for working-class men (Anonymous, 'Mechanics' Institution'; Kelly 1952, 18).[5] Notably, women were also welcomed. Watson estimates that 10–15% of Mechanics' Institute members across Britain were female (Watson 2018, 112). By 1851, there were fifty-five institutes in Scotland and 12,554 members, with many offering lectures in mathematics (Hudson 1851, vi). Dundee Watt Institution had sixty-nine female members by 1850, and the Glasgow Athenaeum had fifty-nine women taking a course of astronomy lectures in October 1847 (73, 82). The origins of the Mechanics' Institutes demonstrate how allowing women into the fringes of scientific pursuits, even with conservative motivations like that of Anderson, could lead to their eventual inclusion in mathematical education.

Another effort to establish opportunities for adult education came from the Glasgow Society for Gaelic Schools, which sought to 'give encouragement and aid by such means and shall appear most effectual for the reading of English, Writing and Arithmetic, in those districts where only Gaelic is known' by opening around nine schools in the Highlands and Islands (Paisley and East Renfrewshire Bible Society 1815, 43). A similar society, the Edinburgh Society for the Support of Gaelic Schools, funded eighty-five such schools (ibid.).[6] In 1816, a visitor to one school at Glencalvie in the Highlands described how 'the mother of the infant is one of the scholars, and such was her ardour to learn that she brought the child and cradle to school' (Hudson 1851, 12). These schools emerged as part of an effort to facilitate the reading of the scriptures in Gaelic-speaking regions, which had been somewhat neglected by the Church; as such, this type of religious education took priority. However, the schools may also have had the incidental benefit of providing lessons in arithmetic to adult women in rural areas.

---

[5] For a recent literature review on Mechanics' Institutes and an assessment of women's participation in them, see (Watson 2018).
[6] For a general account of the ESSGS see (Ritchie 2016).

So far, we have discussed the increasing opportunities for women's mathematical education outside the home, but there were also informal possibilities within the home. The previous quotation from Ferguson highlights another way in which girls might learn arithmetic or mathematics: via their brothers' education. The fictional Eudosia was taught directly by her brother. The English scientific lecturer and instrument maker, Benjamin Martin, used a similar device four years later: Euphrosyne learned 'philosophy' from her brother who was home from university for the vacation (Martin 1772; Millburn 2004). In real-life examples, Joanna Baillie, the Scottish poet, recorded being taught arithmetic by her brother, while Mary Somerville was helped by her brother's tutor to obtain textbooks including Euclid (McMillan 1999, 94–95; Somerville 1874, 52). It is probable that periodicals such as the *Ladies' Diary* provided another route for some women to engage with mathematics. Studies by Shelley Costa, Teri Perl and others have shown how women in England interacted with the *Ladies' Diary*, and it was also read to some extent in Scotland and likely by women (Costa 2002; Perl 1979). So far, however, we have not come across specifically Scottish evidence for such reading, apart from Somerville's recollection of having her curiosity aroused by finding an algebra puzzle in a similar magazine (which was probably not the *Ladies' Diary*). Magazines and books on mathematics appear in the library holdings and many of the private book collections of gentry and professional classes from Dumfriesshire surveyed by Vivienne Dunstan, and are likely to have been available to some young women, though evidence that they were read by either men or women is scanty (Dunstan 2010, 173).

The evident increase in opportunities for women of all social classes across Scotland to gain at least some education in mathematics should be seen in the context of an autodidact culture and rapidly rising literacy and numeracy among both sexes. In the eighteenth century, Lowland Scotland, in particular, achieved one of the highest literacy rates in Europe, with all social classes featuring in library borrowing records and book subscription lists (Rose 2010; Dunstan 2010, 71). Participation rates in mathematics were much lower for both sexes, but women and girls might learn arithmetic or mathematics in parish or private schools or mechanics' institutes, from private tutors or public lecturers. The overwhelming public justification, in Scotland as in the rest of Britain, for providing women with such an education was to fit them as better wives and mothers. Less formally, women might learn within the household context from other family members or through books and magazines, and here there is slightly more evidence of women's own agency in seeking mathematical knowledge in a spirit of pure intellectual curiosity. In either case, the nexus of women, mathematics, and family becomes significant. We discuss this in the next section.

**Family**

The extent to which women encountered mathematics during the eighteenth and nineteenth centuries was often wholly dependent on their domestic situation. During the early eighteenth century, women in both rural and urban areas mostly existed within a household-family model that consisted of parents, children, extended family members and household servants, where 'production and family life were inseparably intertwined', according to Louise Tilly and Joan Scott (Tilly and Scott 1989, 12). Although gender was still important to family hierarchy, the division of labour would occur according to need, rather than to strongly demarcated roles attributed to each sex. However, as the century progressed and industrialization took hold, salaried work outside the home became the burden of men, and 'housework, including family maintenance and childcare, was seen as the appropriate activity for women' (Quataert 1985, 1124). As Danielle Van Den Heuvel notes in her study of early modern Dutch food markets, the 'male-breadwinner model replaced the family economy' (Heuvel 2008, 218). Coinciding with the onset of Enlightenment ideas, the family economy subsided, and new discourses that determined a woman's role within the family emerged. Rosalind Carr and Katie Barclay, in discussing gender in eighteenth-century Scotland, emphasize the 'significance of the idea of women as "helpmeets" for men in shaping women's roles' (Barclay and Carr 2018, 179). Enlightenment philosophers saw virtuous women as being 'those who understood their social role was relational and their purpose to facilitate men' (180). This section will identify how women acquired and utilized computational, and occasionally higher mathematical, knowledge within the family, particularly through accounting, and discuss how this ideal of the 'helpmeet' shaped the circumstances in which such knowledge became acceptable. Women's experiences were also determined by other immediacies like social class, income, and the presence of a husband. Thus, this section will also examine how independent business-owning and widowed women practiced mathematics.

Many women practised accounting and book-keeping as part of a family unit, or could be left in charge of accounts out of necessity in the absence of their husbands. Rosalind Russell has carried out a great deal of research into the lives of business-owning Scottish women during the Enlightenment period. Russell notes that a 'number of Scottish women were self-employed, either setting up their own business or carrying on their husbands' after they had been widowed' (Russell 1988, 4). Several women are known to have run taverns in Edinburgh, and widowed women occupied many kinds of businesses from printing to plumbing, all of which likely involved a level of numeracy to carry out the accounting aspect of the business successfully (5–6). Robert McNair (1703–1779) and Jean Holmes shared equal partnership in their grocery business 'Robert McNair, Jean Holmes & Company', and Holmes herself was trusted with managing the accounts of the firm (Maver, 'Jean Holmes'). Although little is known about her background, it appears that during the 1790s a widowed woman named Mary Brown was 'Glasgow's principal cotton broker', a position which undoubtedly involved calculation (Maver, 'Textile Industry'). In 1815, Susan Carnegie from Montrose created one of the first Savings Banks in Britain for the lower classes, controlling all its financial affairs independently (Ewan et al. 2007, 68). These women joined a growing class of merchants and traders during this period, stimulated by the financial support gained by the Union of 1707 and subsequent Scottish involvement with the British Empire (Carr 2014, 1). Heuvel's observation of the Dutch food markets that the 'processes of commercialization generally benefited independent female entrepreneurship' would appear to be in evidence here too (Heuvel 2008, 217). The independence that these women earned as widows and business-owners created greater utility for mathematical ability.

The use of practical numeracy in the home was not confined to women of lower social classes. According to Elizabeth Ford, Lady Grisell Baillie's account books are meticulously 'detailed records of household management in early eighteenth-century Scotland, especially regarding servants, travel, and the education of her daughters, Rachel (1696–1772) and Grisell (1692–1759)' (Ford 2018, 6). Landlords' absenteeism often led to female aristocrats taking on accountancy roles, like that of Lady Margaret, Countess of Dumfries. She managed the finances of her husband's estate and stated in 1779 that she distributed a wage 'amongst about 30 people in proportion to the number of days they worked at it', indicating that she understood basic arithmetic (Russell 1988, 6). Another aristocrat, Margaret Stewart Calderwood (1715–1774), occupied an estate at Polton, near Edinburgh. She studied mathematics under Professor Colin Maclaurin,[7] and Dorothy McMillan notes that 'the way in which she increased the rental of the family estates would suggest that she profited from the study' (McMillan 1999, 43). Willielma Campbell, Viscountess of Glenorchy (1741–1786), was widowed in 1771 and left with £1,000. She was an 'astute financial manager', founding churches, schools, and an Edinburgh chapel (Russell 1988, 34). The frequency of absentee husbands amongst Scottish women was exacerbated by the Jacobite uprising of 1745 (18). This unsuccessful attempt to restore the Catholic Stuart family to the throne led to the execution and exile of prominent Jacobite leaders, leaving their bereft wives to handle their affairs. One such woman was Lady Isabella Strange, who was left in charge of her exiled husband's finances (7). Her brother described her keen financial management in a letter to her husband, stating that.

> she gives us an account of the present state of the stocks. No Stock-jobber could have done it more distinctly than she does. She talks of per cents., annuities, brokerage, etc., as learnedly as any of the sons of Jonathan (ibid.)

Thus, the political context of eighteenth-century Scotland provided an opportunity and an incentive for women to acquire an understanding of practical mathematics. Although mathematical ability was generally considered to be an inappropriate and unnecessary trait for eighteenth-century Scottish women, in the circumstance of an absent husband, Enlightenment social expectations demanded that a wife obtain the requisite mathematical knowledge to serve her husband's needs as well as possible; to fulfil the role of 'helpmeet', as Barclay and Carr term it.

The demand for household arithmetical skills amongst Scottish women during this period is suggested by the existence of *Institutes of Practical Arithmetic* by William Gordon (1720 or 1721–1793), which was designed as a 'Textbook for Young Ladies' (Mempham 1988, 153; Gordon 1793). Gordon was a Scottish accountant who founded Mercantile Academies in both Glasgow and Edinburgh, and wrote several successful texts on accountancy and bookkeeping (Mempham 1988, 159). Writing in 1793, Gordon emphasized that for a woman who 'pursues any line of business, the knowledge of figures, and accounts, are indispensably necessary' (Gordon 1793, ix). He described several scenarios in which it might be important for women to acquire numerical ability:

> The wife of a merchant, tradesman, or mechanic, as her attention to her husband's business may be frequently necessary, whilst he is absent, or otherways engaged, and still more should she be let [*sic*] a widow, and the whole weight of the business devolve upon herself. (ibid.)

---

[7]Maclaurin is remembered today for the pseudonymous 'Maclaurin series' and was best known in the eighteenth century for his defence of Newton and attempt to introduce classical Greek rigour to calculus in his *Treatise on Fluxions* (1742).

This explicitly framed the desire for women to acquire accounting skills and understanding in terms of her husband's needs and success. The circumstances in which women should learn arithmetic are outlined as those that would serve to benefit her husband or family's enterprise, rather than her own initiative: to serve as her husband's 'helpmeet'. The book covers all basic elements of arithmetic: the number system, units, operators, percentages, ratios, and multiplication and division of numbers (9–10). It is oriented towards practical purposes, so sections include tables of currency conversion, information on the purchases of stocks, insurance premiums, percentages, and ratios to work out proportional profits in partnerships (78). The 'Practical Questions' at the end of each section involve domestic situations, like measuring the length of a carpet, dividing yarn, and requesting cuts of linen, which may indicate the practical applications for which young women were expected to use these skills (10). Although Gordon's sentiments were not particularly progressive in outlook, the existence of the book is an indication that women were indeed expected to learn accountancy for the purposes of assisting their husband's prosperity.

Another notable example of women encountering mathematics through a family environment is that of Edward Sang (1805–1890) and his daughters, Flora (1838–1925) and Jane (1834–1878). Sang was a mathematician from Kirkcaldy, Fife, who published many works on mathematical and engineering problems throughout the nineteenth century (Craik 2003, 47). He is best known for his books of logarithm tables, a task which took forty years and culminated in forty-seven volumes of work (ibid.). He was assisted in this project by Flora and Jane, and a table indicating the 'exact extent of the assistance rendered by Dr Sang's daughters' cites Jane as having compiled five volumes, Flora sixteen volumes, and their father twenty-six volumes (56). Whilst all the 28-place logarithms were computed by Edward, twenty-one volumes out of the thirty-two of 15-place logarithms were calculated by Jane and Flora (ibid.). Sang published books on *Elementary Arithmetic* (1856) and *Higher Arithmetic* (1857), detailing methods of mental calculation that Alexander Craik deems 'very likely' he had taught his daughters, who received much of their education from their parents (51). Although Sang's publications lie just outside the period of interest to this paper (1730–1850), his work is still significant in that it implies that his daughters received a high level of mathematical education during the first half of the nineteenth century, and possessed sufficient mathematical understanding to be able to compute 15-place logarithms to a high degree of accuracy suitable for publication. A statement by his daughter, Flora, indicates that Edward bore no embarrassment at having employed the help of his daughters, and 'did not see why he should not acknowledge the assistance we had given him' (56).

The case of the Sang family exemplifies the idea that it became acceptable and necessary for women to acquire mathematical knowledge if their domestic situation demanded it. Flora and Jane Sang learnt mathematics to contribute to the family enterprise, which in this case was mathematical computation. Although this type of structure was declining by 1850, the circumstances of the family resemble other mathematical family enterprises elsewhere in Britain and Europe, which drew women into the same form of work as their male counterparts. Mary Edwards (1750–1815) was an English woman who carried out almost all of her husband's astronomical computation work for the British *Nautical Almanac*, and trained her daughter Eliza to continue it (Croarken 2003, 9–13). In a similar situation in Germany, Maria Winkelmann-Kirch worked with her husband Gottfried on astronomical observations and computation, and educated their daughters Christina, Maria and Theodora to do the same (Mommertz 2005, 160). Monika Mommertz points out how the 'family members, including the children, as they grew older, were all involved in the working process', a distinctive feature of such a family enterprise (166). As academia and computational work became more formally institutionalized, such practices of working within the home were less prevalent by 1850. However, the Sang family are a later example of how such an enterprise functioned, and demonstrate that such situations enabling women to practice mathematics existed in Scotland, as they did elsewhere.

Family provided many women during the period with the opportunity or need to practice mathematics, at least at the level of keeping household accounts. Although women were largely confined to a domestic sphere that became ever more circumscribed, as industrialization drove an increasing divide between business and home, some businesses, such as the Sangs', remained essentially home-based, enabling discreet assistance by women at a higher mathematical level. During the period, the perception of such assistance shifted from that of contributions to a family economy to that of fulfilling the role of 'help-meet' to the male head. Only in widowhood or, earlier in the period, in the (sometimes enforced) absence of a husband, did women tend to head the enterprise in their own right. Thus, the rationale for educating women in mathematics, to make them better wives and mothers, may have seemed fulfilled. However, as discussed in the next section, some women at least, engaged in mathematics for fun or intellectual curiosity.

**Sociability**

What literature there is on women's engagement with higher mathematics prior to the twentieth century tends to focus on a few outstanding individuals (for example, Mary Somerville), portraying them as isolated figures, unrepresentative of other women. Yet the Enlightenment period was notable for its conversational intellectual culture, encouraged by the popularity of the Parisian salon (Goodman 1989, 331). Many educated women, like Somerville, were embedded in conversational networks through which they might hear of, or discuss, mathematical ideas. These social spaces were characterized by the encouragement of academic conversation and an 'insistence on mixed-gender sociability' (Prendergast 2015, 2). The period between 1730 and 1850 saw a great deal of change in the nature of such conversation. James Secord has demonstrated how, from the 1850s onwards, scientific conversation came to be seen as 'talking shop' due to its association with the trades of the lower commercial classes (Secord 2007, 145). Prior to this, however, 'polite science', including certain aspects of mathematics, formed an appropriate topic of discussion for many women (131). The scientific and literary elite, particularly Somerville and her friends, frequently engaged in polite conversation on mathematics. In her study of progressive science educationalists, Ruth Watts discussed the existence of 'networks in science reaching to Scotland, England and continental Europe' through the formation of friendships in such social spaces (Watts 2007, 97). This section will expand upon this concept, arguing that this network directly facilitated the spread of mathematical knowledge amongst elite women, and can be used to identify Scottish women who may have held mathematical or scientific knowledge.

One woman who engaged in such discussion was Mary Somerville, undoubtedly the most prolific British female mathematician of this period. It is important to understand how her influence could have extended to introduce mathematical ideas to other women in her social circles. After Mary's marriage to her second husband, William, she associated with various friends in salons and informal gatherings (Brück 2009, 73). One friend of the couple was Thomas Young (1777–1829), a physicist, and his wife Eliza, who reportedly was a 'clever woman who shared her husband's interest in mathematics' (ibid.). Mary also associated with Mary Kater (née Reeve, 1784–1829), the wife of Henry Kater, who designed a pendulum which could measure the strength of gravity in various locations (Anonymous, 'Kater's Pendulum'). She was described as being of 'great assistance to her husband'; evidently, she possessed some level of mathematical knowledge (Brück 2009, 74). Both Kater and Somerville were friends with Maria Edgeworth, an Irish author of a scientific book series for children (MacDonald 1977, 94). Edgeworth attended salons in Edinburgh hosted by Elizabeth Hamilton, a novelist and progressive educationalist raised in Scotland (Prendergast 2015, 46). Hamilton's salons were also frequented by Scottish poet Joanna Baillie, who, completing this intellectual 'circle', was also close friends with both Somerville and Edgeworth (76). The presence of scientific women at all these occasions likely engendered mathematical discussion amongst them. On the publication of her translation of Laplace's *Mécanique Céleste*, Somerville gifted copies of the introduction, the 'Preliminary Dissertation', to her female literary friends (Brück 2009, 78). Edgeworth noted in May 1832 that 'you have enlarged my conception of the sublimity of the universe, beyond any ideas I had ever before been enabled to form', and Baillie wrote around the same time that she had 'done more to remove the light estimation in which the capacity of women is too often held, than all that has been accomplished by the whole sisterhood of poetical damsels and novel-writing authors' (78). Although the 'Preliminary Dissertation' contained very basic mathematical concepts, it is interesting that it was read and understood by two women with mostly literary backgrounds, and could potentially have instigated mathematical conversation amongst them (Somerville 1874, 206).

In one mathematical anecdote, Joanna Baillie describes, at the age of twenty-seven, overhearing a friend of ours, a Mathematician, talking one day about squaring the Circle as a discovery which had often been attempted but never found out, and naturally supposing that it must be the discovering a Circle and a Square, exactly the same size, I very simply set my wits to work to find it out (McMillan 1999, 94–95).

To this end, she borrowed 'from my friend Miss Fordyce, now Lady Bentham, an old copy of Euclid', and enlisted her help (ibid.). This anecdote is remarkable for several reasons. Joanna was introduced to this mathematical problem after overhearing a friend discussing it, which would indicate that these networks of friendship did indeed involve mathematical discussion, despite Joanna's literary background. Her description of puzzling over the problem is evidence of a Scottish woman engaging in mathematical problem-solving purely for her own pleasure. She also notes that she borrowed a copy of Euclid's *Elements* from her friend. This was Miss Mary Sophia Fordyce, daughter of Scottish scientist George Fordyce (Pease-Watkin 2002, 3). She married the English engineer and naval architect, Samuel Bentham, in 1796, and it appears that she shared his 'interests in science and engineering', which might explain her possession of a copy of Euclid (ibid.). This extraordinary exchange between two Scottish women, neither of whom had any notable background in mathematics, is exactly the type of mathematical engagement that understanding these informal circles of friendship can uncover.

Not only did this conversational culture allow women to access mathematical knowledge, but it has been argued that the sociable nature of elite women like Mary Somerville encouraged the acceptance of exceptional female ability in mathematics by society. Secord has noted how Somerville may have feared being seen as a 'bizarre specimen' by polite society, so she avoided demonstrations of her extensive knowledge and was instead praised for her 'amiable and gentle' conversation (Secord 2007, 143). Elizabeth Patterson argues that her

> superior domestic management, the warm-hearted sociability and the thoughtful, loving care of children displayed by Mrs Somerville were to return unexpected dividends in the next decades: the 'womanliness' which she exhibited in these ways made her more acceptable as a 'scientific lady' to both the public at large and her scientific peers (Patterson 1969, 319).

Thus, it was not necessarily a more liberal societal attitude towards a woman's role in mathematics, nor even the demonstration that Mary was an extraordinary mathematician, that led to her acceptance in nineteenth-century scientific society. Rather, it was the idea that she could develop mathematical ability in addition to retaining her feminine roles as a 'good wife and mother and a gracious companion' (ibid.). She was permitted success insofar as she remained domestic and feminine. Her sociability also helped her in terms of practical access to scientific resources, through the mediation of her husband, William. Brigitte Stenhouse has demonstrated how William enabled her research by taking on the 'roles of Somerville's chaperone, secretary, representative, or even literary agent' in various situations (Stenhouse 2021, 10). It was through William's connections in Edinburgh that the couple obtained introductions to elite scientific society in London (9). Their partnership allowed her to 'engage productively and meaningfully in the scientific and mathematical communities of which they formed a significant part', irrespective of her gender (16). This climate of conversation between male intellectuals, their wives, and their friends formed a distinctive intellectual culture in eighteenth and nineteenth-century Britain that allowed a woman like Mary Somerville to engage with contemporary mathematical research, despite operating outside of formal institutions.

The formation of such social circles provides connections that enable identification of further Scottish women who are likely to have possessed mathematical knowledge. Jane Davy (née Kerr, Apreece) was a Scottish woman who established an intellectual salon in 1809 in Heriot Row, Edinburgh (Holmes 2008, 243). She turned down a marriage proposal from the Chair of Mathematics and later Professor of Natural Philosophy at the University of Edinburgh, Professor John Playfair (ibid.). Instead, she married Sir Humphry Davy in 1812, a Cornish chemist who lectured for the Royal Institution ('Humphry Davy (1778–1829)', Royal Institution). English writer Sydney Smith described her as the 'first woman who had ever fallen a victim of algebra and been geometrically led from the paths of virtue', and Sir Joseph Banks, an English naturalist, declared 'she has fallen in love with Science and marries him in order to obtain a footing in the Academic Groves' (Lamont-Brown 2007, 73; Holmes 2008, 243). Humphry Davy, and later Jane, was a close personal friend of Maria Edgeworth, and one of Edgeworth's brothers had studied under John Playfair at Edinburgh (Brück 2009, 57). It appears that John Playfair had two sisters, Margaret and Barbara Playfair, who ran a ladies' school in Edinburgh from at least 1806 to 1819, where it is likely that some mathematics was taught (Davies 1872, 53; Anonymous, 'Testament of John Playfair'; 'Collection Guide: Anderson Family Papers M0051, Online Archive of California'). Lady Lucy Clementina Drummond was 'taken to a school kept by Miss Playfair, sister to the professor of that name', just after the birth of her brother, and recalled 'often seeing Sir Humphry Davy, the renowned chemist, alongside Mrs. Apreece, to whom he was subsequently married' (Davies 1872, 53). Although little is known about Jane Davy's educational background, or the curriculum of the Playfairs' school, it seems probable that both involved mathematics to some extent. This school evidently also provided an opportunity for young women to associate with elite intellectuals in mathematics and science, who presumably were present due to the connection with John Playfair. Playfair's protégé, the physicist David Brewster, courted his wife Juliet Macpherson while she was staying with the Misses Playfair in the winter of 1809–10, and the artist Jane Waldie (later Watts) recalled meeting John Playfair frequently while she was a pupil at the school (Gordon 1870, 71–72; Anonymous, 'Mrs Watts', 124). John Playfair's niece, Isabella Best (née Playfair), was also educated at this school, and it appears to have informed and inspired her sufficiently to prompt her to establish her own school for young women when she emigrated to Nova Scotia, Canada (Davies 1995, 14). The Mrs and Misses Best School offered the 'conventional academic and decorative subjects of the day', which likely included some form of mathematics (ibid.).

Another woman who was acquainted with these intellectual circles emanating from Scotland was Jane Marcet (1769–1858). Alongside her husband, a London chemist, she became firm friends with both Mary Somerville and Maria Edgeworth, and attended chemistry and physics lectures at the Royal Institution in London under Sir Humphry Davy and Thomas Young, whose classes were open to women at Young's insistence (Brück 2009, 62–63). It was through one of William's connections in Edinburgh, Leonard Horner, whom he likely met through the Royal Society of Edinburgh, that the Somervilles were introduced to the Marcets (Stenhouse 2021, 4). Inspired by the Royal Institution lectures, Marcet wrote two science books for women, *Conversations on Chemistry: Intended More Especially for the Female Sex* (1805) and later *Conversations on Natural Philosophy* (1826). *Conversations on Chemistry* was highly popular, and is said to have inspired Michael Faraday to establish the Royal Institution's Christmas Lectures for children in 1827, a practice which continues to this day (ibid.). *Conversations on Natural Philosophy* inspired yet another Scottish connection, as after its publication Marcet met with David Brewster and his wife Juliet (64). Marcet had discussed Brewster's theory of colour in her book, and Scottish writer Elizabeth Grant of Rothiemurchus describes entertaining the 'clever authoress' Marcet and the Brewsters at her home in the Highlands in the autumn of 1820, a situation which may have instigated mathematical conversation (Grant 1898). The conversational nature of the book itself is reflective of the idea that Marcet obtained her knowledge through conversation with such intellectual friends (Watts 2007, 92). James Secord notes a 'striking continuity between the kind of topics that could be discussed at such intimate gatherings and the conversational literature of the 1820s and 1830s' (Secord 2007, 140). It is inevitable that in Marcet's understanding of topics including gravity, laws of motion and mechanics, optics, the wave theory of light, and theories of colour, she must have possessed some appreciation for the mathematics involved (Brück 2009, 64). However, what was considered to be significant about the *Conversations* series at the time was the 'accuracy in which all these matters are conveyed in simple language, without the aid of mathematics or even a single formula', and it was intended to be a 'mathematics-free exposition' to science (ibid.). Although the subtitle was removed in later editions of the book, the first was 'intended more especially for the female sex'. Even one of the most well-educated female scientific elite appears to have authored her books under the assumption that women would not be able to grasp the mathematical concepts involved, though whether she believed they lacked the ability, or merely the necessary education, is unclear. Ruth Watts argues that Marcet's book came under a style of 'popular science' that by the 1830s had become 'identified with women and thus were kept in the margins of "real science"' (Watts 2007, 97). This might indicate that although this conversational culture allowed women to enter mathematical discussion and acquire knowledge in a practical sense, it did not change ingrained societal attitudes towards women in mathematics, which still relegated their understanding of science to an exclusively 'female' realm.

**Conclusion**

The extent to which women received a mathematical education, practiced mathematics, and engaged in mathematical research has been overlooked by existing literature on eighteenth and nineteenth-century Scotland. As in other countries, the period between 1730 and 1850 saw the formation of many routes by which women might engage in computational or mathematical activity, both formally and informally, though the balance between informal and formal changed.

Although societal attitudes and expense remained as barriers to women's mathematical education, analysis of newspaper archives has shown that there was a marked proliferation of ladies' schools in Scotland during this period that offered classes in arithmetic, and additionally several in mathematics. Boarding schools in London gave upper-class Scottish women an opportunity to be educated in a metropolis that may have been more sympathetic to their mathematical education. It could be posited that such schools allowed women to cultivate the type of social connections referred to in the 'Sociability' section of this paper, as the accounts of the Playfair sisters' school imply the presence of several elite scientific men and women. Public lectures, mechanics institutes, and Gaelic schools offered other routes into adult education. More work remains to ascertain the level of mathematical education that these girls received, which could be achieved by more intensive studies of school curricula, pupil exercise books, or perhaps personal diaries and commonplace books located in archives: although scholars of reading history have noted a scarcity of useful information in the latter sources.[8]

Women's education remained framed in terms of their usefulness as mothers and homemakers, and thus it is difficult to attribute the slightly greater availability of mathematical education to any progressive Enlightenment change in societal attitudes towards women. One common thread, echoed by both Ferguson and the Scottish

---

[8]For instance, see (Dunstan 2010).

Institution, appears to be that advocacy of women's mathematical education was premised on the idea that mathematics was a rational way for them to occupy their time. Rather than turn to vulgar, unfeminine pursuits like cards, mathematics could help enrich a woman's role as a mother. This perspective still located a woman's place as being within the home, but rather than deem her to be biologically incapable of scientific thought, mathematics became an acceptable and enjoyable pastime.

That women largely practiced mathematics through family and household roles, mainly in the form of accounting and book-keeping, appears to validate the common educational justification, but highlights the barriers erected by cultural expectations. Mandated by circumstances in which a husband might be absent, such as widowhood or separation, women embraced greater economic independence and it became necessary for them to acquire numerical skills. Women in commercial classes might share in the accounting of a business, whilst upper-class women were known to have managed financial affairs for family estates. Families that still relied on a household enterprise structure, like the Sang family, might train young women in the trade, which in this case was mathematics. The Sang family bears a resemblance to other European families who worked in astronomical computation. The presence of accounting textbooks targeted towards women in the late eighteenth century is also indicative of a growing expectation in Scottish society during this period that women should obtain arithmetical abilities, particularly to aid a husband's business venture. Barclay and Carr's analysis of Enlightenment discourse around 'helpmeets' is a useful way in which to understand the conditions under which it became acceptable for women to practice mathematical skills.

Individuals who had a level of mathematical knowledge that went beyond basic arithmetic can be identified in the elite social circles emanating from Scotland, created by the intellectual culture of the Enlightenment. Mixed-gender spaces for polite conversation allowed the wives of male intellectuals to engage in scientific and mathematical discussion, and literary women like Joanna Baillie to be introduced to mathematical ideas in an informal environment. This culture of sociability, in both Scotland and England, allowed Somerville to access contemporary scientific research, and her engagement with it gave her a reputation as an amiable person despite her unusual mathematical ability. Other female scientific women also participated in social circles that involved Scotland, such as Jane Davy and Jane Marcet.

Scottish women's engagement with mathematics through private or social reading of books and magazines has been difficult to assess. While we know that mathematical works were available to many women, discovering reading records will be difficult. Both Vivienne Dunstan and Jill Dye remark on the scarcity of records by readers of either sex of what they were reading, and that book borrowing or purchasing does not equate to readership (Dunstan 2010; Dye 2018). Nor has the Scottish periodical literature been systematically surveyed for its mathematical content. The *Scots Magazine* contained sporadic mathematical puzzles in the 1750s, and routinely included notices and reviews of mathematical books, not differentiated from other forms of literature.[9] This is a potential area for further exploration and would greatly enhance our understanding of how Scottish women engaged with mathematics for leisure if future searches were to prove successful.

---

[9]For example, see 'A Solution of J.C.'s question [xv. 489. 603], and a new question', *Scots Magazine*, 4 February 1754, 75; 'Scottish Literary Intelligence', 1 December 1809, 923.

It is hoped that this study has made a sizeable contribution to a much-neglected field in the history of mathematics in Scotland. Evidently, there was a substantial number of women with mathematical and arithmetical skills in Scotland during this period, perhaps far more than has previously been assumed. There appears to have been some resistance to the idea that women were biologically incapable of mathematical study, though this view persisted, and the subject was still generally considered to be unnecessary and inappropriate for women. Those who did prove this assumption to be incorrect, such as Mary Somerville, were still excluded from formal institutions and viewed as exceptional cases. Certain elements of Scottish and Enlightenment culture did, however, promote women's practice of mathematics at some level; namely, high literacy rates, absent male landowners, and the formation of mixed-gender social spaces through which women could acquire mathematical knowledge and engage in intellectual discussion. Analysis of these social circles remains somewhat under-explored. It is only through understanding all these media that we can obtain a full picture of the extent of mathematical ability amongst women in Scotland during this period.

**Disclosure statement**

No potential conflict of interest was reported by the author(s).

**ORCID**

*Amie Morrison*  http://orcid.org/0000-0001-6678-7203
*Isobel Falconer*  http://orcid.org/0000-0002-7076-9136